\documentclass[11pt,a4paper,twoside]{amsart}

\usepackage[english]{babel}
\usepackage[latin1]{inputenc}
\usepackage{tikz}
\usepackage{pgf}
\usepackage{amssymb}
\usepackage{amsmath}
\usepackage{amsthm}
\usepackage{enumerate}
\usepackage[textsize=tiny]{todonotes} 

\usepackage{mathabx}  
\usepackage{accents}     

\theoremstyle{plain}
\newtheorem{thm}{Theorem}

\newtheorem{defin}[thm]{Definition}
\theoremstyle{remark}

\newtheorem{rmk}[thm]{Remark}

\renewcommand{\epsilon}{\varepsilon}
\newcommand{\ggamma}{\Gamma}

\DeclareMathOperator{\supp}{supp}
\DeclareMathOperator{\sgn}{sgn}

\newcommand{\la}{\langle}
\newcommand{\ra}{\rangle}

\newcommand{\abs}[2][{}]{\lvert{#2}\rvert_{#1}}    

\newcommand{\normsymb}{\|}

\newcommand{\bignormsymb}[1]{#1|\!#1|}
\newcommand{\norm}[2]{\normsymb{#1}\normsymb_{#2}}  
\newcommand{\Bignorm}[2]{\bignormsymb{\Bigl}{#1}\bignormsymb{\Bigr}_{#2}}

\newcommand{\RR}{\mathbb{R}} 
\newcommand{\NN}{\mathbb{N}} 
\newcommand{\ZZ}{\mathbb{Z}} 
\newcommand{\TT}{\mathbb{T}} 

\DeclareMathOperator{\dd}    {d\!}  

\title[Sharp condition for null-controllability on $\RR^d$]{Sharp geometric condition for null-controllability of the heat equation on $\RR^d$ and consistent estimates on the control cost}

\author{Michela Egidi, Ivan Veseli\'c}
\address{Technische Universit\"at Dortmund, 44221 Dortmund, Germany}

\begin{document}
\date{To appear in  \emph{Archiv der Mathematik} with DOI :10.1007/s00013-018-1185-x.}

\maketitle
\begin{abstract}
In this note we study the control problem for the heat equation on $\RR^d$, $d\geq 1$, with control set $\omega\subset\RR^d$.
We provide a necessary and sufficient condition (called $(\gamma, a)$-\emph{thickness}) on $\omega$ such that the heat equation is null-controllable in any positive time.
We give an estimate of the control cost with explicit dependency on the characteristic geometric parameters of the control set.
Finally, we derive a control cost estimate for the heat equation on cubes with periodic, Dirichlet, or Neumann boundary conditions,
where the control sets are again assumed to be thick. We show that the control cost estimate is consistent with the $\RR^d$ case.
\end{abstract}
 \setcounter{tocdepth}{1}


%

\section{Model, definitions, and results}
It is well-known that any open control set $\omega$ is sufficient to achieve null-controllability
of the heat equation on any bounded domain in $\RR^d$ in any positive time.
This means that for any initial datum and time interval the solution of the equation can be driven to zero by a control function supported on $\omega$.
Actually, it has been recently proven in \cite{ApraizEWZ-14,EscauriazaMZ-15} that it is sufficient for $\omega$
to be measurable with positive Lebesgue measure.
If the domain of heat conduction is $\RR^d$ itself, certainly neither of these statements is true.
Hence, from the \emph{qualitative} point of view the requirements on the control set in the bounded and unbounded case
are quite different.
A sufficient condition on control sets for the heat equation on the whole of $\RR^d$ has been recently given
in \cite{LeRousseauM-16}.
In the present paper, we identify a less stringent geometric property on $\omega$ which ensures null-controllability on $\RR^d$.
In fact, we show that it is necessary as well.
This property is called $(\gamma, a)$-\emph{thickness}, where the parameters
$\gamma$ and $ a$ encode the geometry of $\omega$, see below for a precise definition.

Whenever one considers the heat flow (or other transport or conduction problems)
on unbounded domains, this is a mathematical idealization of a physical problem which takes place in a bounded domain.
The idealization is justified in situations where e.g.,~the support of the initial data is known to be far away from the boundary of the finite domain, or if the domain is very large compared to characteristic length scales of the problem.
Therefore it is natural to ask whether the problem on $\RR^d$ and the corresponding one on the bounded domain show the same type of behaviour.
In the case at hand, the characteristic length scales of the problem are specified
by the geometric parameters $\gamma$ and $a$ of the $(\gamma, a)$-thick control set $\omega$.
We establish a control cost estimate which is valid if the domain is $\RR^d$ or any sufficiently large cube, say, of size $L$.
The estimate depends only on $\gamma$ and $a$  (and the control time) but not on the size of the domain.
So, from the \emph{quantitative} point of view the two physical models (large bounded domain versus unbounded domain)
can be reconciled if the quantity of interest is an estimate on the control cost.

One could ask whether it is natural to express the control cost in terms of the parameters $\gamma$ and $a$.
In view of the fact that any feasible control set $\omega$ in  $\RR^d$ is necessarily $(\gamma, a)$-thick for some
$\gamma$ and $a$, one sees that this choice is indeed natural.

The paper is structured as follows:
In the remainder of this section we present crucial definitions, like the one of null-controllability and the thickness of a set,
our main results and a discussion of further questions on the approximability of the heat control problem on large domains by the heat control problem on $\RR^d$.
In Sections \ref{s:thick} and \ref{s:BMPS} we recall fundamental information about thick sets and an abstract observability estimate suitable for our situation, respectively.
Section \ref{s:necessary} contains the proof of the necessity of thickness of the observability set in $\RR^d$ and
Section \ref{s:cost} the control cost estimates, both for $\RR^d$ and for cubes, which use some technical calculations which are deferred to the Appendix.

\subsection{The heat equation, null-controllability, and control cost}

Let $T>0$, $\Omega$ an open region in $\RR^d$, and $S$ any subset of $\RR^d$. We consider the controlled heat equation with control set $\omega=S\cap\Omega$ as
defined below.

\begin{equation}\label{eq:heat_equation}
\left\{\begin{array}{ll}
\partial_t u(t,\cdot) -\Delta u(t,\cdot) = \chi_{\omega}(\cdot)v(t,\cdot) & \qquad\text{ on } [0,T]\times\Omega\\
u(0,\cdot)=u_0(\cdot)\in L^2(\Omega) & \qquad\text{ in } \Omega,
\end{array}
\right.
\end{equation}
where $\Delta$ denotes the Laplacian on $\Omega$
and $\chi_{\omega}$ the characteristic function of $\omega$.
If $\Omega \neq \RR^d$, we will always specify below which type of boundary conditions are imposed on $\partial\Omega$.
The function $v\in L^2([0,T]\times \Omega)$ is called \emph{control function}.

\begin{defin} The system \eqref{eq:heat_equation} is null-controllable in time $T>0$ if for every initial datum $u_0\in L^2(\Omega)$
there exists a control function $v\colon [0,T]\times\Omega\to \RR$, $v\in L^2([0,T]\times\Omega)$ such that the solution of \eqref{eq:heat_equation} satisfies $u(T, \cdot) \equiv 0$.

The quantity
\begin{multline}
C_T:=\sup_{\norm{u_0}{L^2(\Omega)}=1}\inf\{\norm{v}{L^2([0,T]\times \omega)}\;\vert\; \text{the  solution $u$ of \eqref{eq:heat_equation} with r.h.s. } \\
\chi_\omega\cdot v \text{ satisfies } u(T,\cdot)\equiv 0\}
\end{multline}
is called control cost.
\end{defin}
The main purpose of this work is to identify a necessary and sufficient condition on the set $\omega$
such that the heat equation \eqref{eq:heat_equation} on $\RR^d$ is null-controllable from  $\omega$ in any time $T>0$.
The key geometric property is the following.

\begin{defin}
Let $\gamma\in (0,1]$ and let $a=(a_1,\ldots, a_d)\in\RR^d_+$.
Set $A:= [0,a_1]\times \cdots\times[0,a_d]\subset\RR^d$.
A set $S\subset\RR^d$ is called  $(\gamma,a)$-thick if it is measurable and
\[
\abs{S\cap \left(x+A\right)  }\geq \gamma \prod_{j=1}^d a_j,
\]
 for all $x\in \RR^d$. Here $\abs{\cdot}$ stands for the Lebesgue measure on $\RR^d$.
A set $S\subset\RR^d$ is called  thick if there exist $\gamma\in (0,1]$ and $a\in\RR^d_+$ such that $S$ is  $(\gamma,a)$-thick.
\end{defin}
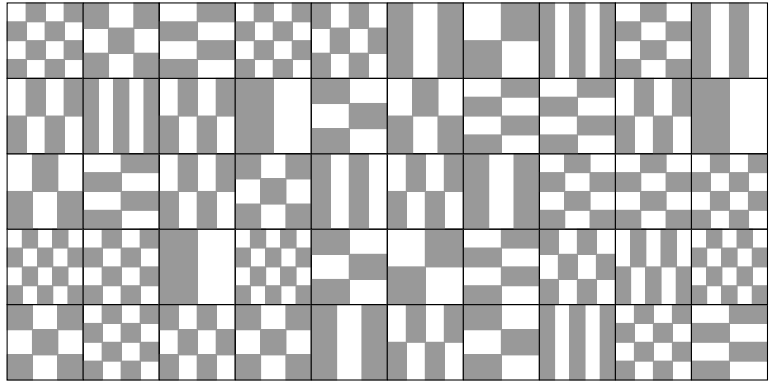
\begin{figure}[ht]\centering
 \begin{tikzpicture}
  \pgfmathsetmacro{\X}{9};
  \pgfmathsetmacro{\Y}{4};
  \pgfmathsetmacro{\e}{5};
  \pgfmathsetmacro{\d}{4};
  \foreach \y in {0,...,\Y}{
	\foreach \x in {0,...,\X}{
		\pgfmathsetmacro{\A}{random(2,\e)}; 
		\pgfmathsetmacro{\a}{\A-1};
                \pgfmathsetmacro{\B}{random(1,\d)};  
		\pgfmathsetmacro{\b}{\B-1};
		\pgfmathsetmacro{\Lo}{1/\A}; 
		\pgfmathsetmacro{\Lv}{1/\B}; 
		\foreach \i in {0,..., \a}{
			  \pgfmathsetmacro{\v}{\x +(\i *\Lo)};
			  \pgfmathsetmacro{\V}{\x +(\i +1)*\Lo};
			  \foreach \j in {0,...,\b}{
				\pgfmathsetmacro{\W}{\y + (\j +1)*\Lv};
				\pgfmathsetmacro{\w}{\y +(\j *\Lv)};
						\pgfmathsetmacro{\test}{\j +\i};
						\ifthenelse{\isodd{\test}}{
						}{\filldraw[black!40] (\v,\w) rectangle (\V, \W);}
			}
		}
		\draw (\x,\y) rectangle (\x+1, \y+1);
	}
  }
 \end{tikzpicture}
\caption{Example of a section of size $10 \times 5$ of a $(\gamma, a)$-thick set $S$ (marked in gray) in  $\RR^2$ where $a=(2,2)$ and $\gamma=1/8$.}
\end{figure}
In particular, the set $S$ needs not to be open.

\subsection{Results on null-controllability on $\RR^d$ and on cubes}
\begin{thm}\label{thm:null-control-full-space}%
 Let $T>0$ and consider problem \eqref{eq:heat_equation} on $\RR^d$ with control set $\omega=S$. The following statements are equivalent.
\begin{itemize}
 \item [(i)] $S\subset\RR^d$ is a thick set ;
 \item [(ii)] the heat equation on $\RR^d$ is null-controllable in time $T>0$.
\end{itemize}

Furthermore, if $S$ is $(\gamma, a)$-thick in $\RR^d$ the control cost satisfies
\begin{equation}\label{eq:control-cost-full-space}
 C_T\leq C_1^{1/2}\exp\left(\frac{C_1}{2T}\right),
\quad \text{ where } \quad
C_1= \left(\frac{K^d}{\gamma}\right)^{K(d+\norm{a}{1})} 
\end{equation}
for $\norm{a}{1}=a_1+\ldots+ a_d$, $K=12\sqrt{2} K_1$, and $K_1$ the universal constant of
Remark \ref{remark:constants}.
\end{thm}

\begin{thm}\label{thm:control-heat-torus}
Let $T,L>0$, $\Omega=\TT^d_L:= (0, 2\pi L)^d$ and $S$ be a $(\gamma, a)$-thick set in $\RR^d$ with $a_j\leq 2\pi L$ for all $j=1,\ldots, d$.
Let $\Delta$ be the Laplacian on $\Omega$ with periodic, Dirichlet, or Neumann boundary conditions on $\partial\Omega$ and $\omega=S\cap\TT^d_L$. Then, the heat equation \eqref{eq:heat_equation} on $\TT^d_L$
with control set $\omega$ is null-controllable in time $T>0$ and the control cost satisfies
\begin{equation}\label{eq:control-cost-cubes}
 C_T\leq C_2^{1/2}\exp\left(\frac{C_2}{2T}\right),
\quad \text{ where } \quad
 C_2=  \left(\frac{K^d}{\gamma}\right)^{K(d+\norm{a}{1})},
\end{equation}
for $\norm{a}{1}= a_1+\ldots +a_d$, $K=24\sqrt{2}K_1$, and $K_1$ the universal constant of Remark \ref{remark:constants}.
\end{thm}

A crucial ingredient in the proof of Theorem \ref{thm:null-control-full-space} is the so-called
Logvinenko-Sereda Theorem. We use a quantitative version of it due to \cite{Kovrijkine-thesis}.
In a previous paper \cite{EgidiV-16} we have derived an analogous result for finite cubes, used in
Theorem \ref{thm:control-heat-torus}.
Independently of us \cite{WangWZZ} obtained the equivalence
(i) $\Leftrightarrow$ (ii) in Theorem \ref{thm:null-control-full-space}.

If one assumes that the thick set $S$  satisfies additionally that there exist $G, \delta>0$ such that for each $k \in (G\ZZ)^d$
the intersection of $S$ and the cube $k+[0,G]^d$ contains a $\delta$-ball then it is possible to derive control cost estimates similar to those in
Theorem \ref{thm:null-control-full-space} and 
Theorem \ref{thm:control-heat-torus} 
also for the null-control problem for Schr\"odinger semigroups, see \cite{NakicTTV-18} 
and \cite{NakicTTV-18-preprint}, respectively for details.

\subsection{Discussion, questions, and further goals concerning the approximation of control functions on $\RR^d$}

Our results about consistent
estimates on the control cost for the limiting heat control problem on $\RR^d$ and the
approximating problems on cubes $\Lambda_L$ sets the stage for many questions and further research.
Note in particular that our consistent control cost estimates are uniform in the
size $L$ of the cubical domain $\Lambda_L$.

If $v_L$ is the null control function for the problem on
$\Lambda_L$ with minimal $L^2([0,T]\times \Lambda_L)$-norm and correspondingly
$v$ the null control function for the problem on
$\RR^d$ with minimal $L^2([0,T]\times \RR^d)$-norm,
we expect that $v_L$ converges to $v$ as $L \to \infty$.
Even more interesting is a quantitative bound on $v_L-v$.
A promising approach to implement this is to study the variational functionals
associated with Lions' Hilbert Uniqueness Method, which characterize the controls with minimal norm
and give a quite explicit description using the adjoint system.
\smallskip

Alternatively, one can ask what is the $L^2(\RR^d)$-norm of $w(T, \cdot)$,
if $w$ is the solution of
\begin{equation}
\left\{\begin{array}{ll}
\partial_t w(t,\cdot) -\Delta w(t,\cdot) = \chi_{S}(\cdot)v_L(t,\cdot) & \qquad\text{ on } [0,T]\times\RR^d\\
w(0,\cdot)=u_0(\cdot)\in L^2(\Lambda_L) &
\end{array}
\right.
\end{equation}
where $v_L$ is the (minimum norm) $T$-null control function of
\begin{equation}
\left\{\begin{array}{ll}
\partial_t u_L(t,\cdot) -\Delta u_L(t,\cdot) = \chi_{S_L}(\cdot)v_L(t,\cdot)
& \qquad\text{ on } [0,T]\times\Lambda_L\\
u_L(0,\cdot)=u_0(\cdot)\in L^2(\Lambda_L) &
\end{array}
\right.
\end{equation}
which ensures $u_L(T, \cdot )\equiv 0$. The estimate on $w(T, \cdot)$
will certainly depend on the scale $L$, or more precisely on the
location of the support of the initial condition $u_0$ within the cube $\Lambda_L$.
\smallskip

This and more will be discussed in a sequel project.
\bigskip

Since we prove uniform estimates on the control cost in the scale $L$
it is possible to  extract weakly convergent subsequences from $v_L$.
This leads to questions concerning the
properties of weak accumulation points.
Let us start with a functional analytic observation.

\begin{rmk}
Let $X$ and $Y$ be separable Hilbert spaces and assume
that $H_L\colon X \to Y$ for each $L \in \NN$ are bounded linear operators.
Then the following statements hold true. \smallskip

(i) \quad
\emph{If $c_H:=\sup_{L\in\NN} \|H_L\|< \infty$ then there is a subsequence $K \subset \NN$
and a bounded linear $H:X \to Y$ such that for all $z\in Y$ and $x\in X$
\[
\lim_{K \ni L \to \infty} \la z, H_L x\ra =\la z, H x\ra,
\]
i.e. $\lim_{K \ni L \to \infty} H_L =H$ weakly.}\smallskip

To see this, choose a dense countable subset $D$ of $X$.
For each $x \in D$ the sequence
$H_{L}x $ has a weakly convergent subsequence to some  $y \in Y$. By choosing a diagonal subsequence one
can achieve that along this subsequence the above weak convergence holds for all $x\in  D$ simultaneously.
The map $D \ni x \to y \in Y$ can be shown to be a linear and bounded operator. Consequently, it can be extended to all of $X$,
still as a linear and bounded operator $H:X \to Y$.
Finally, it follows that along the chosen diagonal subsequence,
for \emph{every} $x\in X$ one has $H_{L_k}x \to H x$ weakly in $Y$. \smallskip

(ii) \quad
Now we discuss the (almost) converse of (i).\smallskip

\emph{Assume there is a subsequence $K \subset \NN$
and a bounded linear $H:X \to Y$
such that $\lim_{K \ni L \to \infty} H_L =H$ weakly, i.e.
for all $z\in Y$ and $x\in X$
\[
\lim_{K \ni L \to \infty} \la z, H_L x\ra =\la z, H x\ra.
\]
Then $c_H:=\sup_{L\in K} \|H_L\|< \infty$.}
\smallskip

It is obvious that we can have only uniform boundedness along $K$ and not along $\NN$.
Apart from this the statements form  an equivalence.
\smallskip

To prove (ii), let $K = (L_k)_k$ and define for each $k \in \NN$,
$
T_k: Y\times X \to \RR$ by $T_k(z,x)= \la z, H_{L_k}x\ra
$. Then $T_k$ is bi-linear and bounded,
\[
\|T_k\|:= \sup_{z\in Y, \|z\|=1} \sup_{x\in X, \|x\|=1} |T_k(z, x)| \leq \|H_{L_k} \|,
\]
and by the assumed weak convergence there exists some $k_0(z)\in \NN$ such that
\[
\sup_k |T_k(z, x)| \leq \max\{|\la z, H  x\ra| +1,|T_1(z, x)|, \dots, |T_{k_0}(z, x)|   \} =: c_{z,x}.
\]
Now the multi-linear uniform boundedness principle implies the desired statement. \smallskip
\end{rmk}

Now we turn to some control theoretic aspects of weak convergence and
consider the heat equation on the cube $\Lambda_L$
\begin{equation}\label{eq:heat-equation-L}
\left\{\begin{array}{ll}
\partial_t u_L(t,\cdot) -\Delta u_L(t,\cdot) = v_L(t,\cdot) & \qquad\text{ on } [0,T]\times\Lambda_L\\
u_L(0,\cdot)=u_0(\cdot)\in L^2(\Lambda_L) & \qquad\text{ in } \Lambda_L,
\end{array}
\right.
\end{equation}
where  we impose Dirichlet boundary conditions on $\partial \Lambda_L$.
Then using the Feynman-Kac and Duhamel formulae, one can prove the following.

\begin{rmk}\label{rm:weak-convergence-of-uL}
Let $u_0\in L^2(\RR^d)$ be supported in $\Lambda_R$ for some $R>0$.
For each $L \in \NN, L>R$,
let $v_L \in L^2([0,T]\times \Lambda_L)\subset L^2([0,T]\times \RR^d)$ be the right hand side of
the heat equation  \eqref{eq:heat-equation-L} on $\Lambda_L$.
Assume that
\begin{enumerate}[(i)]
  \item $v_L$ converges weakly to $v \in L^2([0,T]\times \RR^d)$,
  \item $\sup_{L\in \NN, L>R} \|v_L\| < \infty$.
\end{enumerate}
Then the solutions $u_L$ of \eqref{eq:heat-equation-L} converge weakly in $L^2([0,T]\times \RR^d)$  to the solution $u$ of
\begin{equation}
\left\{\begin{array}{ll}
\partial_t u(t,\cdot) -\Delta u(t,\cdot) = v(t,\cdot) & \qquad\text{ on } [0,T]\times\RR^d\\
u(0,\cdot)=u_0(\cdot)\in L^2(\RR^d) & \qquad\text{ in } \RR^d.
\end{array}
\right.\\[1em]
\end{equation}
\end{rmk}

Note that in general there are many $T$-null controls $v_L$ driving the solution of
\eqref{eq:heat-equation-L} to zero.
(If minimality of the norm is assumed, we have uniqueness of $v_L$, but we are considering
all controls whose norm is bounded by the constant $C_T$ in \eqref{eq:control-cost-full-space}
and \eqref{eq:control-cost-cubes}.)
For each weak accumulation point of such a sequence  $v_L, L\in \NN$, as in
Remark \ref{rm:weak-convergence-of-uL} we obtain a corresponding limiting solution $u$.
A priori the set of such functions $u$ can be very large.\smallskip

Finally, let us in the situation of Remark \ref{rm:weak-convergence-of-uL}
additionally assume that $v_L$ is the $T$-null control with minimal norm and that for all $L$ we have
$u_L(T,\cdot)\equiv 0$ on $\Lambda_L$. Under this condition we believe that the weak limit $u$ satisfies
$u(T,\cdot)\equiv 0$ on $\RR^d$ as well.

\section{Thick sets and their properties}
\label{s:thick}

We recall some results on $(\gamma, a)$-thick sets  $S \subset\RR^d$ which are needed in the proof of Theorems \ref{thm:null-control-full-space} and \ref{thm:control-heat-torus}.

The question on which linear subspaces of $L^2(\RR^d)$ the seminorm $f \mapsto \|f\|_{L^2(S)}$ is actually equivalent to the original norm goes back at least to
Panejah, see \cite{Panejah-61,Panejah-62}. A sharp condition was identified almost simultaneously, but using different proofs, by Logvinenko \& Sereda and  Kacnel'son in
\cite{Logvinenko-Sereda-74,Katsnelson-73}. Kovrijkine \cite{Kovrijkine-01,Kovrijkine-thesis} has turned these qualitative equivalence statements into a quantitative estimate.

\begin{thm}[Logvinenko-Sereda theorem by \cite{Kovrijkine-01}]\label{thm:logvinenko-sereda-full-space}
Let $f\in L^2(\RR^d)$ and assume that $\supp\hat{f}\subset J$, where $J$ is a parallelepiped with sides parallel to coordinate axes and of length $b_1,\ldots,b_d$.
Set $b=(b_1,\ldots,b_d)$.
Let $S$ be a $(\gamma, a)$-thick set. Then,
\begin{equation}\label{eq:kov1}
\norm{f}{L^2(\RR^d)}\leq\left(\frac{K_1^d}{\gamma}\right)^{K_1(a\cdot b+d)}\norm{f}{L^2(S)},
\end{equation}
where $a\cdot b$ stands for the euclidean inner product and $K_1$ is a universal constant which w.l.o.g. may be taken bigger or equal than $\mathrm{e}$.
\end{thm}
%

The above Theorem has been adapted for functions on the torus $\TT^d_L$ in \cite{EgidiV-16}.

\begin{thm}[Logvinenko-Sereda on the torus by \cite{EgidiV-16}]\label{thm:Log-Ser_torus}
	Let $f\in L^2(\TT^d_L)$ and assume that $\supp \hat{f}\subset J$, where $J$ is a parallelepiped
	with sides parallel to the coordinate axes and of length $b_1,\ldots,b_d$. Set $b=(b_1,\ldots,b_n)$.
	Let $S\subset\RR^d$ be a $(\gamma, a)$-thick set with $a_j\leq 2\pi L$ for all $j=1,\ldots, d$. Then,
	\begin{equation}\label{eq:Log-Ser_one_cube}
	\norm{f}{L^2(\TT^d_L)}\leq\left(\frac{K_2^d}{\gamma}\right)^{ K_{2} a\cdot b+\frac{6d+1}{2}}\norm{f}{L^2(\TT_L^d\cap S)},
	\end{equation}
	where $a\cdot b$ stands for the euclidean inner product and $K_2$ is a universal constant which w.l.o.g. may be taken bigger or equal than $\mathrm{e}$.
\end{thm}

\begin{rmk}\label{remark:constants}
The estimates \eqref{eq:kov1} and \eqref{eq:Log-Ser_one_cube} are both of the type $\left(\frac{\tilde{K}^d}{\gamma}\right)^{\tilde{K}(d+a\cdot b)}$,
for some universal constant $\tilde{K}$.
In particular, if we increase the constant $K_1$ in Theorem \ref{thm:logvinenko-sereda-full-space} to satisfy
\[
 K_1\geq \max(K_2, 7/2),
\]
then  \eqref{eq:Log-Ser_one_cube} implies
\[
 \norm{f}{L^2(\TT^d_L)}\leq \left(\frac{K_1^d}{\gamma}\right)^{K_1(a\cdot b+d)}\norm{f}{L^2(\TT_L^d\cap S)}.
\]
\end{rmk}

\section{An abstract observability result of
\cite{BeauchardPS-18}}
\label{s:BMPS}

Null-controllability of the system \eqref{eq:heat_equation} in time $T>0$
is equivalent to the following observability estimate with respect to $\omega$,
see for example \cite[Theorem 2.44]{coron:07}.
\begin{equation}\label{eq:observability-def}
 \exists C>0\   \forall g_0\in L^2(\Omega)\, : \quad \norm{g(T,\cdot)}{L^2(\Omega)}^2\leq C\int_0^T \norm{g(t,\cdot)}{L^2(\omega)}^2\dd t,
\end{equation}
where $g$ is the solution of the adjoint system
\begin{equation}\label{eq:heat_equation-adjoint}
\left\{\begin{array}{ll}
 \partial_t g(t,\cdot) - \Delta g(t,\cdot) = 0 & \qquad\text{ on } [0,T]\times\Omega\\
 g(0,\cdot)=g_0\in L^2(\Omega) & \qquad\text{ in }\Omega.
\end{array}\right.
\end{equation}
Here $\Delta$ has the same boundary conditions as in problem \eqref{eq:heat_equation}.
In addition, the equivalence provides an estimate for the control cost, namely
\[
C_T\leq \sqrt{C},
\]
where $C$ is the observability constant appearing in \eqref{eq:observability-def}.

Hence, we can use a special case of a general observability result
obtained by K.~Beauchard, L.~Miller, and K.~Pravda-Starov
in the recent work \cite[Theorem 2.1]{BeauchardPS-18} to show null-controllability.
\begin{thm}\label{thm:observability-B-PS}
Let $\Omega$ be an open subset of $\RR^d$,
$\omega\subset \Omega$ be measurable, $(\pi_E)_{E\in\NN}$ a family of orthogonal projections in $L^2(\Omega)$ ,
$(e^{t\Delta})_{t\geq 0}$ the contraction semigroup associated to a Laplacian (with self-adjoint boundary conditions) on $L^2(\Omega)$,
and $c_1,t_0>0$ positive constants.
If the spectral inequality
\begin{equation}\label{eq:spectral-inequality-general}
\forall f\in L^2(\Omega), \ \forall E\in\NN\, : \quad \norm{\pi_E f}{L^2(\Omega)}\leq e^{c_1 \sqrt{E}}\norm{\pi_E f}{L^2(\omega)},
\end{equation}
and the dissipation estimate
\begin{equation}\label{eq:dissipation-general}
\forall f\in L^2(\Omega), \ \forall E\in\NN, \ \forall t\in (0,t_0)\, : \quad \norm{(1-\pi_E)(e^{t\Delta}f)}{L^2(\Omega)}\leq e^{- t E}\norm{f}{L^2(\Omega)}
\end{equation}
hold, then there exists a positive constant $C_3>1$ such that the observability estimate
\begin{equation}\label{eq:observability-general}
\forall T>0, \ \forall f\in L^2(\Omega)\, : \quad \norm{e^{T\Delta}f}{L^2(\Omega)}^2\leq C_3\exp\left(\frac{C_3}{T}\right)\int_0^T \norm{e^{t\Delta}f}{L^2(\omega)}^2 \dd t
\end{equation}
holds true.
\end{thm}

We remark that in the original statement in \cite{BeauchardPS-18}
the set $\omega$ was assumed to be open, however the proof works also for measurable subsets.

\section{Necessary condition for null-controllability on  $\RR^d$ in Theorem \ref{thm:null-control-full-space} }
\label{s:necessary}

Due to the duality between observability and null-controllability, it is sufficient to show that
the adjoint system \eqref{eq:heat_equation-adjoint} is not observable in time $T>0$ if $S$ is not a thick set.
The latter condition means that
\[
\forall \gamma >0, a \in \RR_+^d \quad \exists \xi=\xi_{\gamma,a}\in \RR^d \ : \quad  |S \cap (\xi+A) | <\gamma.
\]
 In particular, for any $k \in \NN$ there exists a $\xi_k \in \RR^d$ such that $|S \cap (\xi_k+[0,2k]^d) | <1/k$. Thus for $x_k=\xi_k+(k,\ldots,k) \in \RR^d$ we have

\begin{equation}\label{eq:assumption-contradiction}
\abs{S\cap B(x_k, k)}< 1/k,
\end{equation}
where $B(x_k, k)$ denotes the euclidean ball in $\RR^d$ centred at $x_k$ with radius $k$.

We now construct a sequence of functions which shows that \eqref{eq:observability-def} does not hold.
For this purpose, choose the initial data $g_k(0,x):= \exp\left(-\frac{\norm{x-x_k}{2}^2}{2}\right)$. Then, the solution of \eqref{eq:heat_equation-adjoint}
associated to this initial value is
\begin{equation*}
 g_k(t,x)=\frac{1}{(2t+1)^{d/2}}e^{-\frac{\norm{x-x_k}{2}^2}{2(2t+1)}}.
\end{equation*}

By the change of variable $y=x-x_k$ we have
\begin{align*}
 \int_{\RR^d}\abs{g_k(t,x)}^2 \dd x & = \int_{\RR^d}\frac{1}{(2t+1)^{d}}e^{-\frac{\norm{x-x_k}{2}^2}{(2t+1)}} \dd x
 = \frac{1}{(2t+1)^d}\int_{\RR^d}e^{-\frac{\norm{y}{2}^2}{(2t+1)}} \dd y.
\end{align*}
Thus the $L^2$-norm of $g_k(t,x)$ in $\RR^d$ is independent of $k$.
However, the term $\int_0^T \norm{g_k(t,\cdot)}{L^2(S)}^2 \dd t$ converges to zero as $k\to \infty$, as we show now.
\begin{align}
\nonumber
\int_0^T\int_S \frac{1}{(2t+1)^d} e^{-\frac{\norm{x-x_k}{2}^2}{2t+1}} &\dd x \dd t   =
\int_0^T\int_{S- x_k} \frac{1}{(2t+1)^d} e^{-\frac{\norm{y}{2}^2}{2t+1}} \dd y \dd t\\
\nonumber
 & \leq \int_0^T\int_{S-x_k} e^{-\frac{\norm{y}{2}^2}{2T+1}} \dd y \dd t\\
\nonumber
 & \leq T\int_{B(0, k)\cap (S-x_k)} e^{-\frac{\norm{y}{2}^2}{2T+1}} \dd y  +T \int_{\norm{y}{2}\geq  k} e^{-\frac{\norm{y}{2}^2}{2T+1}} \dd y \\
\nonumber
 & \leq T\abs{B(0,k)\cap(S-x_k)} + T\int_{\norm{y}{2} \geq k} e^{-\frac{\norm{y}{2}^2}{2T+1}} \dd y \\
\label{eq:null-sequence}
 & = T\abs{B(x_k,k)\cap S} + T\int_{\norm{y}{2}>k} e^{-\frac{\norm{y}{2}^2}{2T+1}} \dd y ,
\end{align}
where we used the substitution $y=x-x_k$, the monotonicity of the exponential in $t$, and the estimate
$e^{-\frac{\norm{y}{2}^2}{2T+1}}\leq 1$.
Now, both terms in \eqref{eq:null-sequence} go to zero as $k\to\infty$ and, consequently, so does $\int_0^T \norm{g_k(t,\cdot)}{L^2(S)}^2 \dd t$, which  yields the desired contraposition.

\section{Control cost estimates}
\label{s:cost}
In this section we derive explicit control cost estimates for the heat equation when the observability set is $(\gamma, a)$-thick
for some $\gamma >0$ and $a \in \RR^d_+$. This implies, in particular, that such sets are sufficient for null-controllability of the heat equation.
We first derive the control cost estimate for Theorem \ref{thm:null-control-full-space} and then the one for Theorem  \ref{thm:control-heat-torus}.

\begin{proof}[Control cost estimate for  Theorem \ref{thm:null-control-full-space}]
We assume that $S$ is $(\gamma, a)$-thick and we show that the adjoint system \eqref{eq:heat_equation-adjoint} is
observable in time $T>0$ using Theorem \ref{thm:observability-B-PS}.
To this end, we only need to check its assumptions, namely the dissipation and spectral inequality.

Let $E\in \NN$ and $\pi_E=\chi_{(-\infty,E]}(-\Delta)$, so that $1-\chi_{(-\infty,E]}(-\Delta)=\chi_{(E,+\infty)}(-\Delta)$.
By spectral calculus $e^{2t\Delta}\chi_{[E,+\infty)}(-\Delta) \leq e^{-2tE}\chi_{[E,+\infty)}(-\Delta)$ in the sense of quadratic forms.
This yields the dissipation estimate
\begin{equation}
\label{eq:dissipation-estimate}
\begin{aligned}
\norm{\chi_{(E,+\infty)}(-\Delta)(e^{t\Delta}f)}{L^2(\RR^d)}^2
& =\langle \chi_{(E,+\infty)}(-\Delta) f, \chi_{(E,+\infty)}(-\Delta) e^{2t\Delta}f\rangle_{L^2(\RR^d)}\\
& \leqslant e^{-2tE}\norm{\chi_{(E,+\infty)}(-\Delta)f}{L^2(\RR^d)}^2\\
& \leq e^{-2tE}\norm{f}{L^2(\RR^d)}^2,
\end{aligned}
\end{equation}
which implies \eqref{eq:dissipation-general}.

We set $J=[-\sqrt{E},\sqrt{E}]^d$, $b=(2\sqrt{E},\ldots, 2\sqrt{E})$, and deduce
the spectral inequality for $\pi_E f$ from Theorem \ref{thm:logvinenko-sereda-full-space}.
Note that $ \supp \widehat{\pi_E f } \subset J$, hence we have
\begin{align*}
\norm{\pi_E f}{L^2(\RR^d)} &
\leq \left(\frac{K_1^d}{\gamma}\right)^{K_1(d+2\sqrt{E}\norm{a}{1})}\norm{\pi_E f}{L^2(S)}
\\
&\leq\exp\left(2\sqrt{E}K_1(d+\norm{a}{1})\log\left(\frac{K_1^d}{\gamma}\right)\right) \norm{\pi_E f}{L^2(S)},
\end{align*}
and so inequality \eqref{eq:spectral-inequality-general}
holds with
\begin{equation}\label{eq:lower-bound-c1}
c_1=2K_1(d+\norm{a}{1})\log\left(\frac{K_1^d}{\gamma}\right) \geqslant 2 d^2 \mathrm{e}
\end{equation}
since $K_1 \geqslant \mathrm{e}$.
By Theorem \ref{thm:observability-B-PS}, there exists a constant $C_3>1$ such that for all $T>0$
and for all $f\in L^2(\RR^d)$ we have
\begin{equation*}
\norm{e^{T\Delta}f}{L^2(\RR^d)}^2\leq C_3\exp\left(\frac{C_3}{T}\right)\int_0^T \norm{e^{t\Delta}f}{L^2(S)}^2\; \dd t.
\end{equation*}
More precisely, from the Appendix we infer $C_3=\exp(6\sqrt{2} c_1)$ and use
\eqref{eq:lower-bound-c1} to conclude
\begin{equation*}
C_3
= \left(\frac{K_1^d}{\gamma}\right)^{12\sqrt{2}K_1(d+\norm{a}{1})}
\leq \left(\frac{K^d}{\gamma}\right)^{K(d+\norm{a}{1})},
\text{where $K=12\sqrt{2}K_1$}.
\end{equation*}
This leads to the control cost estimate
$ C_T\leq C_3^{1/2}\exp\left(\frac{C_3}{2T}\right)$.
\end{proof}


The proof of Theorem \ref{thm:control-heat-torus} proceeds along the same lines.
The main difference is that we have to use Theorem \ref{thm:Log-Ser_torus}
instead of Theorem \ref{thm:logvinenko-sereda-full-space}
and use an expansion in terms of the eigenbasis of the Laplacian with periodic, Dirichlet, or Neumann boundary conditions.

\begin{proof}[Control cost estimate on the cube and proof of Theorem \ref{thm:control-heat-torus}]
We only need to check the dissipation estimate and spectral inequality of Theorem \ref{thm:observability-B-PS} for $\pi_E=\chi_{(-\infty, E]}(-\Delta)$, $E\in\NN$.
Since the Laplacian on $\TT^d_L$ with periodic, Dirichlet, or Neumann boundary conditions is self-adjoint,
the dissipation estimate \eqref{eq:dissipation-estimate} holds verbatim.

To derive the spectral inequality we first observe the following.
Any $L^2$-functions on $\TT^d_L$ is a linear combination of eigenfunctions of the periodic, Dirichlet or Neumann Laplacian on $\TT^d_L$, as they all form an orthonormal
basis for $L^2(\TT^d_L)$. Such eigenfunctions are
\begin{align*}
 \phi_k^P(x)=  & \left(2\pi L\right)^{-d/2} e^{i \frac{k}{L}\cdot x}\quad k\in\ZZ^d\\
  \phi_k^D(x)= & \left(\pi L\right)^{-d/2}\prod_{j=1}^d\sin\left(\frac{k_j x_j}{2 L}\right) \quad k\in\NN^d,\\
  \phi_k^N(x)= & \left(\pi L\right)^{-d/2}\prod_{j=1}^d\cos\left(\frac{k_j x_j}{2 L}\right)\quad k\in\NN_0^d,
\end{align*}
for periodic, Dirichlet and Neumann boundary conditions, respectively, corresponding to eigenvalues
%
\[
\lambda_k=\frac{\norm{k}{2}^2}{L^2},\quad \norm{k}{2}^2=\sum_{j=1}^d \abs{k_j}^2, \quad k\in\ZZ^d
\]
for periodic boundary conditions and
\[
\lambda_k=\frac{\norm{k}{2}^2}{(2L)^2},\quad \norm{k}{2}^2=\sum_{j=1}^d \abs{k_j}^2
\]
with $k\in\NN^d$ for Dirichlet boundary conditions, and $k\in\NN_0^d$ for Neumann boundary conditions.

Set now $J=[-\sqrt{E},\sqrt{E}]^d$ and consider the case of periodic boundary conditions. Let $f=\sum_{k \in \ZZ^d}\alpha_k^{P}\phi_k^{P}$, then
\begin{equation*}
\pi_E f(x)=
\sum_{k \in \ZZ^d, \sqrt{\lambda_k} \leq \sqrt{E}}\alpha_k^{P} \left(\frac{1}{2\pi L}\right)^{d/2}e^{i \frac{k}{L}\cdot x}
=\left(\frac{1}{2\pi L}\right)^{d/2}\sum_{\frac{k}{L}\in J\cap\left(\frac{1}{L}\ZZ\right)^d}a_k^P e^{i\frac {k}{L}\cdot x},
\end{equation*}
where
\[
a_k^P=\left\{
\begin{array}{ll}
\alpha_k^P & \frac{k}{L}\in B_{\sqrt{E}}(0)\cap\left(\frac{1}{L}\ZZ\right)^d\\
0 & \text{ otherwise }.
\end{array}\right.
\]
Consequently, by Theorem \ref{thm:Log-Ser_torus} we calculate
\begin{align*}
 \norm{\pi_E f}{L^2(\TT^d_L)}^2 &
= \Bignorm{\left(\frac{1}{2\pi L}\right)^{d/2}\sum_{\frac{k}{L}\in J\cap\left(\frac{1}{L}\ZZ\right)^d}a^P_k e^{i \frac{k}{L} } }{L^2(\TT^d_L)}^2\\
& \hspace{-1cm} \leq \left(\frac{K_{2}^d}{\gamma}\right)^{ 4K_{2}\sqrt{E}\norm{a}{1}+6d+1}\Bignorm{\left(\frac{1}{2\pi L}\right)^{d/2}\sum_{\frac{k}{L}\in J\cap\left(\frac{1}{L}\ZZ\right)^d} a_k^P e^{i \frac{k}{L}} }{L^2(\TT^d_L\cap S)}^2\\
& \hspace{-1cm} \leq \left(\frac{K_{1}^d}{\gamma}\right)^{ 4K_{1}\sqrt{E}(\norm{a}{1}+d)}
\norm{\pi_E f}{L^2(\TT^d_L\cap S)}^2,
 \end{align*}
 for $K_1$ as in Remark \ref{remark:constants}. Then, the spectral inequality is fulfilled with
 \[
  c_1 = 2K_1( \norm{a}{1}+d)   \log  \left(\frac{K_{1}^d}{\gamma}\right) \geq 2d^2\mathrm{e},
 \]
since $K_1\geq \mathrm{e}$.

We now consider the case of Dirichlet and Neumann boundary conditions. Let $f^D=\sum_{k \in \NN^d}\alpha_k^{D}\phi_k^{D}$ and $f^N=\sum_{k \in \NN_0^d}\alpha_k^{N}\phi_k^{N}$, then
\begin{align*}
& \pi_E f^D(x)=\sum_{k \in \NN^d, \sqrt{\lambda_k} \leq \sqrt{E}}\alpha_k^{D}\left(\frac{1}{\pi L}\right)^{d/2}\prod_{j=1}^d\sin\left(\frac{k_j x_j}{2 L}\right),\\
& \pi_E f^N(x) = \sum_{k \in \NN_0^d, \sqrt{\lambda_k} \leq \sqrt{E}}\alpha_k^{N}\left(\frac{1}{\pi L}\right)^{d/2}\prod_{j=1}^d\cos\left(\frac{k_j x_j}{2 L}\right).
\end{align*}
We observe that $\widehat{\pi_E f^D}$ and $\widehat{\pi_E f^N}$ have no compact support.
However, when $\pi_E f^D$ and $\pi_E f^N$ are extended to the following $L^2$-functions on $\TT^d_{2L}$
\begin{align*}
 &\widetilde{\pi_E f^D}(x)=\sum_{k \in \NN^d, \sqrt{\lambda_k} \leq \sqrt{E}}\alpha_k^{D}\left(\frac{1}{\pi L}\right)^{d/2}\prod_{j=1}^d\sin\left(\frac{k_j x_j}{2 L}\right),\qquad x\in\TT^d_{2L},\\
 &\widetilde{\pi_E f^N}(x)=\sum_{k \in \NN_0^d, \sqrt{\lambda_k} \leq \sqrt{E}}\alpha_k^{N}\left(\frac{1}{\pi L}\right)^{d/2}\prod_{j=1}^d\cos\left(\frac{k_j x_j}{2 L}\right),\qquad x\in\TT^d_{2L},
\end{align*}
their respective Fourier Series on $\TT^d_{2L}$ have Fourier coefficients concentrated in $J$.

To see this, embed the ball $B_{\sqrt{E}}(0)$ of radius $\sqrt{E}$ centred at zero into the $d$-dimensional cube $J$ and use the relations $\sin x=(2i)^{-1}(e^{ix}-e^{-ix})$ and $\cos x= 2^{-1}(e^{ix}+ e^{-ix})$ to obtain
\begin{align*}
& \widetilde{\pi_E f^D}(x)
=\left(\frac{1}{\pi L}\right)^{d/2}\left(\frac{1}{2i}\right)^d\sum_{\frac{\tilde{k}}{2L}\in J\cap\left(\frac{1}{2L}\ZZ\right)^d}a^D_{\tilde{k}} e^{i \frac{\tilde{k}}{2L}\cdot x},\\	
& \widetilde{\pi_E f^N}(x)=\left(\frac{1}{\pi L}\right)^{d/2}\left(\frac{1}{2}\right)^d\sum_{\frac{\tilde{k}}{2L}\in J\cap\left(\frac{1}{2L}\ZZ\right)^d}a^N_{\tilde{k}} e^{i \frac{\tilde{k}}{2L}\cdot x},
\end{align*}
where
\[
 a^D_{\tilde{k}}=\left\{
  \begin{array}{ll}
  \alpha_{\abs{\tilde{k}}}^D\prod_{j=1}^d\sgn(\tilde{k_j}) & \qquad\frac{\tilde{k}}{2L}\in B_{\sqrt{E}}(0)\cap\left(\frac{1}{2L}\ZZ\right)^d,\\
  0 & \qquad\text{otherwise},
 \end{array}
  \right.
\]
and
\[
a^N_{\tilde{k}}=\left\{
  \begin{array}{ll}
  \alpha_{\abs{\tilde{k}}}^N & \qquad\frac{\tilde{k}}{2L}\in B_{\sqrt{E}}(0)\cap\left(\frac{1}{2L}\ZZ\right)^d,\\
  0 & \qquad\text{otherwise}.
 \end{array}
  \right.
\]
for $\abs{\tilde{k}}=(\abs{\tilde{k_1}},\ldots, \abs{\tilde{k_d}})\in\NN_0^d$. This means that the Fourier coefficients of the extensions are all contained in $J$.
We define
\begin{align*}
  S^{(0)} &= S\cap [0,2\pi L]^d \\
  S^{(1)} &= S^{(0)} \cup \{(-x_1,x_2,\ldots,x_d)\mid (x_1,x_2,\ldots,x_d)\in   S^{(0)} \} \\
  S^{(2)} &= S^{(1)} \cup \{(x_1,-x_2,x_3,\ldots,x_d)\mid (x_1,x_2,\ldots,x_d)\in   S^{(1)} \}\\
  \vdots  &= \vdots \\
  S^{(d)} &= S^{(d-1)} \cup \{(x_1,\ldots,x_{d-1},-x_d)\mid (x_1,x_2,\ldots,x_d)\in   S^{(d-1)} \}
\end{align*}
 and extend $S^{(d)}$ periodically to $ \tilde{S}=  \bigcup_{\kappa \in (4\pi L\ZZ)^d} \left(\kappa +S^{(d)}\right) $.
 We claim that $\tilde{S}$ is a $(\gamma/2^d, 2a)$-thick set.

Let $A:= [0,a_1]\times \cdots\times[0,a_d]$, $\tilde A:= [0,2a_1]\times \cdots\times[0,2a_d]$, and $x\in\RR^d$ arbitrary.
If $(x+\tilde A)^{\circ}$ is disjoint to the lattice $\ggamma= (2\pi L \ZZ)^d$, then it is contained in a periodicity cell
$\kappa + [0,2\pi L]^d$ with $\kappa \in \ggamma$. This is still true for $(x+A)^{\circ}$. Since the Lebesgue measure does not charge boundaries of parallelepipeds and
is reflection invariant, we have
\[
 \abs{(x+\tilde A)^{\circ}\cap \tilde S }= \abs{(x+\tilde A)\cap \tilde S }\geq  \abs{(x+A)\cap \tilde S } \geq
\gamma \prod_{j=1}^d a_j = \left(\frac{\gamma}{2^d}\right) \prod_{j=1}^d (2a_j).
\]
If $(x+\tilde A)^{\circ}$ contains a lattice point $\nu \in \ggamma$, then there is a periodicity cell
$\kappa + [0,2\pi L]^d$ with $\kappa \in \ggamma, \Vert \nu -\kappa \Vert \leq 1$ such that
 $(x+\tilde A)^{\circ}\cap (\kappa+ [0,2\pi L]^d)$ is a parallelepiped whose $j$th side has lenght at least $a_j$.
In particular, it contains $y+A$ for some $y \in \RR^d$. As above we have
\[
 \abs{(x+\tilde A)^{\circ}\cap \tilde S } \geq  \abs{(y+A)\cap \tilde S } \geq
\gamma \prod_{j=1}^d a_j = \left(\frac{\gamma}{2^d}\right) \prod_{j=1}^d (2a_j),
\]
which concludes the proof of the claim.

We now set $\Gamma_1 =\{0, 2\pi L\}^d$ and, to unify notation, $*\in\{D,N\}$. 
By reflection symmetry of $\widetilde{\pi_E f^*}$ and $\tilde S$, and by change of variable we have
$\norm{\widetilde{\pi_E f^*}}{L^2(\tilde S \cap (\kappa +\TT^d_{L})) }^2
=\norm{\widetilde{\pi_E f^*}}{L^2(\tilde S \cap (\nu +\TT^d_{L})) }^2$
and
$\norm{\widetilde{\pi_E f^*}}{L^2(\kappa +\TT^d_{L}) }^2
=\norm{\widetilde{\pi_E f^*}}{L^2(\nu +\TT^d_{L})}^2$ if $\kappa, \nu\in \ggamma$ are neighbouring lattice points.
Consequently
\begin{align*}
\norm{\widetilde{\pi_E f^*}}{L^2(\TT^d_{2L})}^2 &
=\sum_{\kappa\in \Gamma_1}\norm{\widetilde{\pi_E f^*}}{L^2(\kappa+\TT_{L}^d)}^2=2^d\norm{\pi_E f^*}{L^2(\TT^d_L)}^2,\\
\norm{\widetilde{\pi_E f^*}}{L^2(\tilde{S}\cap \TT^d_{2L})}^2 &
=\sum_{\kappa\in \Gamma_1}\norm{\widetilde{\pi_E f^*}}{L^2(\tilde{S} \cap (\kappa+\TT_{L}^d))}^2
=2^d\norm{\pi_E f^*}{L^2(S\cap \TT^d_L)}^2.
\end{align*}
Now, Theorem \ref{thm:Log-Ser_torus} applied to $\widetilde{\pi_E f^*}$ and $\tilde{S}$ yields
\begin{align*}
 \norm{\pi_E f^*}{L^2(\TT^d_{L})}^2 & = 2^{-d}\norm{\widetilde{\pi_E f^*}}{L^2(\TT^d_{2L})}^2\\
 & \leq 2^{-d}\left(\frac{(2 K_2)^d}{\gamma}\right)^{ 8 K_{2}\sqrt{E} \norm{a}{1}+6d+1}\norm{\widetilde{\pi_E f^*}}{L^2(\tilde{S}\cap \TT^d_{2L})}^2\\
 & = \left(\frac{(2 K_2)^d}{\gamma}\right)^{ 8 K_{2}\sqrt{E} \norm{a}{1}+6d+1} \norm{\pi_E f^*}{L^2(S\cap \TT^d_{L})}^2\\
 & = \left(\frac{(2 K_1)^d}{\gamma}\right)^{ 8 K_{1}\sqrt{E} (\norm{a}{1}+d)} \norm{\pi_E f^*}{L^2(S\cap \TT^d_{L})}^2,
\end{align*}
for $K_1$ as in Remark \ref{remark:constants}, and the spectral inequality \eqref{eq:spectral-inequality-general} is satisfied for
\[
 c_1= 4 K_1( \norm{a}{1}+d)\log \left(\frac{(2 K_{1})^d}{\gamma}\right) \geq 4d^2\mathrm{e}.
\]

From the above estimate and the Appendix we can now conclude that the observability estimate \eqref{eq:observability-general} is satisfied with constant
\[
C_3= \exp(6\sqrt{2}c_1) \leq \left(\frac{K^d}{\gamma}\right)^{K(\norm{a}{1}+d)} ,\quad K=24\sqrt{2}K_1,
\]
 which also gives
$
 C_T\leq C_3^{1/2}\exp\left(\frac{C_3}{2T}\right).
$
\end{proof}

\appendix

\section{Observability constant $C_3$}
We here show how to compute the constant $C_3$ in Theorem \ref{thm:observability-B-PS}.
We will use the notation from \cite[Appendix 8.3]{BeauchardPS-18}
(apart from the numbering of the constants $C_j$), where \cite[Theorem 2.1]{BeauchardPS-18}
is proven.
Recall that $c_1\geq 2d^2 \mathrm{e}$ and
define\footnote{In the notation of \cite{BeauchardPS-18} we choose $q:=1/2,m=1, a=1/2,b=1, c_2=1$
thus $\gamma=(3 c_1)^2 2^5$ but these parameters do not appear in the sequel. }
\[
M:=8(3 c_1)^2
\]
We are looking for some  $\tau_0$ such that for all $0<\tau <\tau_0$
the followings are fulfilled
\begin{align}
&\frac{(3 c_1)^2 2^5}{\tau^2}>1,\label{eq:1}\\
&\frac{1}{\tau}\exp\left(-\frac{2^3 3 c_1^2}{\tau}\right)\leq \frac{1}{4},\label{eq:2}\\
&\frac{1}{\tau}\exp\left(\frac{2^4 (3 c_1)^2}{\tau}\right)\geq 1.\label{eq:3}
\end{align}
Set $\tau_0:=2^{5/2}3 c_1$.

Note that eq.~\eqref{eq:1} is equivalent to $\tau< \tau_0$.

Eq.~\eqref{eq:2} is fulfilled for any $0<\tau\leq 2^3 3 c_1^2=:\tau_1 $.
In fact, the function $(0, \infty)\ni \tau \to\frac{1}{\tau}\exp\left(-\frac{2^3 3 c_1^2}{\tau}\right)$
has its unique maximum at $\tau_1$.
Hence, for all $0<\tau\leq \tau_1$ and our choice of $c_1$ we have
\[
\frac{1}{\tau}\exp\left(-\frac{2^3 3 c_1^2}{\tau}\right)\leq \frac{1}{\tau_1}\exp\left(-\frac{2^3 3 c_1^2}{\tau_1}\right)=
\frac{1}{2^3 3 c_1^2 \mathrm{e}}\leq \frac{1}{2^3 3 (2 d^2 \mathrm{e})^2 \mathrm{e}}\leq \frac{1}{96\mathrm{e}^3}< \frac{1}{4}
\]
Eq.~\eqref{eq:3} is fulfilled for all $\tau\leq \tau_0$ since the function is decreasing
on $(0,\infty)$ and for  $\tau\in (0,\tau_0)$ we obtain
\[
\frac{1}{\tau}\exp\left(\frac{2^4 (3 c_1)^2}{\tau}\right)\geq
\frac{1}{2^{5/2} 3 c_1}\exp(2^{3/2}3c_1)
\geq \frac{2^{3} (3 c_1)^2 +1}{2^{5/2} 3 c_1} 
\geq 3 c_1 \cdot\sqrt{2}\geq 1,
\]
where we used $\exp(x)\geq 1 +\frac{x^2}{2}$ and $c_1\geq 2d^2 \mathrm{e}>1$.

Consequently, Eq.~\eqref{eq:1}, \eqref{eq:2}, and \eqref{eq:3} are all satisfied
for $0<\tau<\tau_0= 2^{5/2} 3 c_1 $.
Now, from \cite[Appendix 8.3]{BeauchardPS-18}
we extract
\[
 C_3=\max(C_4, C_5)
\]
where $C_4= 2 M = 144 c_1^2$ and $C_5=\exp\left(\frac{2 C_4}{2\tau_0}\right)=\exp(6\sqrt{2} c_1)$. We conclude that
\[
 C_3=\max\left(144 c_1^2, \exp(6\sqrt{2} c_1)\right)=\exp(6\sqrt{2} c_1)
\]
since $\exp(6\sqrt{2} c_1) \geq 144 c_1^2 $ for all $c_1 \geq 0.3086$ and in our case $c_1\geq 2d^2 \mathrm{e}$.

\subsection*{Acknowledgements}
This work has been partially supported by the Deut\-sche Forschungsgemeinschaft under grant no.~VE 253/7-1
\emph{Multiscale version of the Logvinenko-Sereda Theorem}.
Comments on earlier versions of the manuscript and helpful discussions
with I.~Naki\'c, A.~Seelmann, M.~T\"aufer, M.~Tautenhahn
and K.~Veseli\'c are gratefully acknowledged.
We also thank the anonymous referee for stimulating comments on the submitted manuscript.

\bibliography{bibliography}
\bibliographystyle{alpha}
\end{document}